\documentclass[conference]{IEEEtran}
\IEEEoverridecommandlockouts
\usepackage{cite}
\usepackage{amsmath,amssymb,amsfonts}
\usepackage{algorithmic}
\usepackage{graphicx}
\usepackage{textcomp}
\usepackage{xcolor}
\usepackage{epstopdf}
\usepackage{amsthm}
\def\BibTeX{{\rm B\kern-.05em{\sc i\kern-.025em b}\kern-.08em
		T\kern-.1667em\lower.7ex\hbox{E}\kern-.125emX}}
\def\caly{{\mathcal Y}}

\def\calr{{\mathcal R}}

\def\calp{{\mathcal P}}

\def\calz{{\mathcal Z}}

\def\L2e{{\cal L}_{2e}}

\def\rea{\mathbb{R}}

\def\adj{\mbox{adj}}

\newtheorem{assumption}{Assumption}
\newtheorem{statement}{Statement}
\newtheorem{proposition}{Proposition}

\begin{document}

\title{Adaptive observer for a LTV system with partially unknown state matrix and  delayed measurements\\
\thanks{This work is supported by the Russian Science Foundation under grant 22-21-00499.}
}

\author{\IEEEauthorblockN{1\textsuperscript{st} Alexey Bobtsov}
	\IEEEauthorblockA{\textit{Department of Control Systems and Robotics} \\
		\textit{ITMO University}\\
		Saint-Petersburg, Russia \\
		bobtsov@mail.ru}
\and
\IEEEauthorblockN{2\textsuperscript{nd} Nikolay Nikolaev}
\IEEEauthorblockA{\textit{Department of Control Systems and Robotics} \\
	\textit{ITMO University}\\
	Saint-Petersburg, Russia  \\
	nikona@yandex.ru}
\and 
\IEEEauthorblockN{3\textsuperscript{rd} Olga Slita}
\IEEEauthorblockA{\textit{Department of Control Systems and Robotics} \\
	\textit{ITMO University}\\
	Saint-Petersburg, Russia  \\
	o-slita@yandex.ru}
\and

\IEEEauthorblockN{4\textsuperscript{th} Olga Kozachek}
\IEEEauthorblockA{\textit{Department of Control Systems and Robotics} \\
	\textit{ITMO University}\\
	Saint-Petersburg, Russia  \\
	oakozachek@itmo.ru}
}

\maketitle

\begin{abstract}
Problem of adaptive state observer synthesis for linear time-varying (LTV) system with unknown time-varying parameter and delayed output measurements is considered. State observation problem has attracted the attention of many researchers \cite{b4}. In this paper the results proposed in the \cite{b2}, \cite{b8}, \cite{b9} are developed. It is supposed that the state matrix can be represented as sum of known and unknown parts. Output vector is measured with known constant delay. An adaptive identification algorithm which reconstructs unknown state and unknown time-varying parameter is proposed.
\end{abstract}

\begin{IEEEkeywords}
LTV system, delay, adaptive observer.
\end{IEEEkeywords}

\section{Introduction and Problem Formulation}
The common problem in control applications is that real devices provide measurements with delays. This fact makes the problem of the observer and control laws design for dynamical systems more complicated. This problem has been explored in many papers recently. In case of linear time invariant (LTI) systems, this issue is well known \cite{b1}. However, for LTV systems this problem is still widely open - see the literature review and references in the recent papers  \cite{b4}, \cite{b2}, \cite{b3}. Currently, researchers consider different problems for systems with delayed measurements, for instance, estimation of unknown measurement delay  \cite{b12};  state observers for cases with fixed constant delay \cite{b16} and with time-varying delay \cite{b2}, \cite{b8},\cite{b9} ets. As examples of the technical systems with measurement delays authors consider the problem of biomass regulation \cite{b13}, \cite{b14}, autonomous underwater vehicle tracking problem \cite{b15}, actuator fault estimation  \cite{b17}, ets.

In this paper a linear time-varying system with unknown time-varying parameter is considered. An observer of the unknown parameter is designed under condition that the output signal is measured with a constant delay.

In this paper we consider a LTV system described by the following equation
\begin{align}
	\dot{x}(t)=A_0(t)x(t)+B(t)u(t),
\end{align}
where $x(t) \in \rea^n$ is an unknown state vector, $ A_0 (t)\in {\rea}^{n\times n}$ and $B(t)\in {\rea}^{n}$.
\begin{assumption}
	Let us suppose that state matrix $A_0(t)$ can be rewritten in the following form $A_0(t)=A(t)+\theta(t) I$, where   $ I\in {\rea}^{n\times n}$ is the identity matrix.
\end{assumption}

According to the Assumption 1, we can rewrite the original system in the following form
\begin{align}
	\label{sys}
	\dot{x}(t)&=A(t)x(t)+\theta(t)x(t)+B(t)u(t),\\ 
	y(t)&=C(t) x(\phi(t)),
	\label{out}
\end{align}
where $x(t) \in \rea^n$ is an unknown state vector, $ y(t)\in {\rea}^{n} $ is a vector of measured output variables, $ u(t) \in \rea $ is a known input signal, $ A(t)\in {\rea}^{n\times n},B(t)\in {\rea}^{n} $ and $ C(t) \in {\rea}^{n\times n}$ are known matrices. The entries of $A(t)$, $B(t)$ and $C(t)$ are assumed to be continuous and bounded, $ \phi (t) $ is a continuous known nonnegative function which defines the measurement delay
\begin{equation}
		\phi (t)=t-d,\phi(t)\ge 0, 
\end{equation}
where $ d> 0 $ is a constant delay; $ \theta (t) $ is an unknown time-varying function defined by the equation
\begin{equation} 
	\ddot{\theta }(t)=-{\omega }^{2}\theta(t), 
	\label{thet}
\end{equation}
where $ \omega > 0 $ is a constant unknown parameter.

It is obvious that \eqref{thet} is a model of harmonic signal generator, so it can be rewritten in the following form
\begin{align}
	\label{theta}
	{\theta }(t)={{a }}_{1}\sin {\omega }t+{{a }}_{2}\cos {\omega }t,
\end{align}
where $a_1$ and $a_2$ are corresponding constant parameters.
\begin{assumption}
	We suppose that $C(t)$ is the identity matrix $C(t)={I}_{n \times n}$.
\end{assumption}
	
In this paper we design an adaptive observer of unmeasured state $x(t)$ such that the following inequality holds
\begin{align}
	\lim\limits_{t \to \infty} (x(t)-\hat{x}(t)) = 0,
\end{align}
where $\hat{x}(t)$ is the estimate of the state vector of the system \eqref{sys} and  \eqref{out}.

\section{Main Result}

\subsection{Preliminary transformations}

The solution of the problem of adaptive observer design for LTV system can be obtained in two steps. On the first step the constant parameter $ \omega $ is estimated and on the second step  an observer for $ \theta $ is designed.

Let us consider the system (\ref{sys}), \eqref{out} at the moment $ t – d $. The equation \eqref{sys} can be rewritten in the following way
\begin{align}
	\dot{{x}}_d={A}_d {x}_d+\eta(t){x}_d+{B}_d {u}_d,
	\label{sys_delay}
\end{align}
where
\begin{subequations}
	\begin{align}
		{x}_d&=x(\phi (t))=y, \\
		{A}_d&=A(\phi (t)), \\
		{B}_d&=B(\phi (t)), \\
		{u}_d&=u(\phi (t)), \\
		\eta(t)&=\theta (\phi (t)).
		\label{eta_theta} 
	\end{align}
	\end{subequations}

Thus, the new parameter $ \eta(t) $ is a solution of the equation
\begin{equation}
	\ddot{\eta }=-{\omega }^{2}\eta .
	\label{eta}
\end{equation}

Let us define a new function $V$
\begin{equation}
	V=y^\top y = {{x}}_d^\top{x}_d.
\end{equation}

From (\ref{sys_delay}) we can write the derivative of this function in the following form
\begin{equation}
	\dot{V}={{x}}_d^\top({{A}}_d^\top+{A}_d){x}_d+2{u}_d{{x}}_d^\top{B}_d+2\eta (t){{x}}_d^\top{x}_d,
\end{equation}
or
\begin{equation}
	\dot{V}=\alpha ({x}_d)+2\eta (t){{x}}_d^\top{x}_d,
\end{equation}
where
\begin{equation}
	\alpha ({x}_d):={{x}}_d^\top ({{A}}_d^\top+{A}_d){x}_d+2{u}_d{{x}}_d^\top {B}_d
\end{equation}
is a known function.

Let us define a new variable $\xi$ by the equation
\begin{equation}
	\xi =\ln V.
\end{equation}

Its derivative can be written in the following form
\begin{equation}
	\dot{\xi}=\frac{\dot{V}}{V}=\frac{\alpha ({x}_d)}{V}+2\eta(t) .
\end{equation}

From the previous equation the function $ \eta(t) $ can be found as
\begin{equation}
	\eta(t) =\frac{1}{2}(\dot{\xi }+\beta ({x}_d)),
	\label{eta_expr}
\end{equation}
where
\begin{equation}
	\beta({x}_d)=-\frac{\alpha({x}_d)}{V}.
\end{equation}

\subsection{Estimation of unknown parameter $\omega$}

The equation (\ref{eta}) can be rewritten as
\begin{equation}
	{p}^{2}\eta(t) =-{\omega }^{2}\eta(t) ,
	\label{eta_dif}
\end{equation}
where $ p=d/dt $ is a differential operator.

Let us apply the filter $ \frac{\lambda_1^3}{{(p+\lambda_1)}^{3}} $, where $\lambda_1>0$ to (\ref{eta_dif}). Then we obtain
\begin{equation}
	\label{filter}
	\frac{\lambda_1^3{p}^{2}}{{(p+\lambda_1)}^{3}}\eta(t) =-{\omega }^{2}\frac{\lambda_1^3}{{(p+\lambda_1)}^{3}}\eta(t) .
\end{equation}

Substituting $ \eta(t) $ from (\ref{eta_expr}) into this equation, we receive
\begin{equation}
	\frac{\lambda_1^3 {p}^{2}}{{(p+\lambda_1)}^{3}}\frac{1}{2}(\dot{\xi }+\beta({x}_d))=-{\omega }^{2}\frac{\lambda_1^3}{{(p+\lambda_1)}^{3}}\frac{1}{2}(\dot{\xi }+\beta({x}_d)).
\end{equation}

Previous expression can be rewritten as a linear regression equation of the form
\begin{equation}
	q=\varphi k,
\end{equation}
where
\begin{equation}
	q=\frac{\lambda_1^3{p}^{3}}{{(p+\lambda_1)}^{3}}\xi +\frac{\lambda_1^3{p}^{2}}{{(p+\lambda_1)}^{3}}\beta({x}_d), \\
\end{equation}
\begin{equation}
	\varphi =\frac{\lambda_1^3p}{{(p+\lambda_1)}^{3}}\xi +\frac{\lambda_1^3}{{(p+\lambda_1)}^{3}}\beta({x}_d),\\
\end{equation}
\begin{equation}
	k=-{\omega }^{2}.
\end{equation}

The unknown parameter $ k $ can be estimated using, for instance, standard gradient algorithm, see \cite{b5}, \cite{b6}
\begin{equation}
	\dot{\hat{k}}=\gamma_1 \varphi (q-\varphi \hat{k}),
\end{equation}
where $\gamma_1>0$ is an adaptation gain.

Then, the unknown parameter $ \omega $ can be obtained as
\begin{equation}
	\label{hat_omega}
	\hat{\omega }=\sqrt{\left| \hat{k} \right|}.
\end{equation}

Consider the error
\begin{align}
	\tilde{k}=\hat{k}-k.
\end{align}

Then for derivative of $\tilde{k}$ we have
\begin{align}
	\dot{\tilde{k}}=\dot{\hat{k}}=\gamma_1 \varphi (q - \varphi \tilde{k}) =  \gamma_1 \varphi^2 {k} -\gamma_1 \varphi^2 \hat{k} = - \gamma_1 \varphi^2 \tilde{k}.
\end{align}

The solution for $\tilde{k}$ takes the form
\begin{align}
	\tilde{k}(t) = \tilde{k}_0 e^{-\gamma_1 \int\limits_0^t \varphi^2d\tau} 
\end{align}

Then
\begin{align}
	\hat{k}(t)=k+\tilde{k}_0 e^{-\gamma_1 \int\limits_0^t \varphi^2d\tau}
\end{align}
and for $\hat{\omega}$ we receive
\begin{align}
\hat{\omega}(t)=\sqrt{	k+\tilde{k}_0 e^{-\gamma_1 \int\limits_0^t \varphi^2d\tau}}=\omega+\varepsilon(t),
\end{align}
where $\varepsilon(t)$ is exponentially  decaying term due to $\tilde{k}_0 e^{-\gamma_1 \int\limits_0^t \varphi^2d\tau}$.

\begin{statement} Convergence of $\hat{\theta}$ to $\theta$.
	
Consider equation \eqref{theta}. Let us find estimate of $\theta$ in the form $\hat{\theta }(t)=\hat{a}_1 \sin \hat{\omega}t + \hat{a}_2\cos \hat{\omega }t$.

Consider the error $\tilde{\theta}=\hat{\theta}-\theta$.
\begin{align}
	\nonumber
	\tilde{\theta}&=\hat{\theta}-\theta=\hat{a}_1 \sin \hat{\omega}t + \hat{a}_2\cos \hat{\omega }t - {{a }}_{1}\sin ({\omega }t) \\ 
	\nonumber
	&- {{a }}_{2}\cos ({\omega }t) 
	=(\tilde{a}_1+a_1)\sin(\omega t + \varepsilon t ) \\
	\nonumber
	&+ (\tilde{a}_2+a_2)\cos(\omega t + \varepsilon t ) - a_1 \sin(\omega t) - a_2 \cos(\omega t)\\
	\nonumber
	&= \tilde{a}_1 \sin(\omega t + \varepsilon t ) + \tilde{a}_2 \cos(\omega t + \varepsilon t )  \\ 
	\nonumber
	&+ a_1 (\sin(\omega t + \varepsilon t ) - \sin ({\omega }t)) + a_2 (\cos(\omega t + \varepsilon t ) - \cos ({\omega }t))
\end{align}
where $\tilde{a}_i=\hat{a}_i-a_i$, $i=1,2$.

Since $\sin(\omega t + \varepsilon t )=\sin(\omega t) \cos(\varepsilon t)+\cos(\omega t) \sin(\varepsilon t)$ and  $\cos(\omega t + \varepsilon t )=\cos(\omega t) \cos(\varepsilon t)-\sin(\omega t) \sin(\varepsilon t)$, we have
\begin{align}
\nonumber
\tilde{\theta}(t)&=\tilde{a}_1 \sin(\omega t + \varepsilon t ) + \tilde{a}_2 \cos(\omega t + \varepsilon t) \\
\nonumber
& + a_1(\sin(\omega t) \cos(\varepsilon t)+\cos(\omega t) \sin(\varepsilon t)-\sin(\omega t))\\
\nonumber
& + a_2(\cos(\omega t) \cos(\varepsilon t) - \sin(\omega t) \sin(\varepsilon t)-\cos(\omega t)).
\end{align}
It is easy to show that 
\begin{align}
	\nonumber
	\lim\limits_{t \to \infty} \cos (\varepsilon t) = 1,\\
	\nonumber
	\lim\limits_{t \to \infty} \sin (\varepsilon t) = 0,\\
	\nonumber
	\lim\limits_{t \to \infty} \tilde{a}_i(t) = 0,
\end{align}
so finaly we have
\begin{align}
\nonumber
\lim\limits_{t \to \infty} \tilde{\theta}(t) = a_1(\sin(\omega t)-\sin(\omega t)) + a_2(\cos(\omega t)-\cos(\omega t)).
\end{align}
\end{statement}

\subsection{Time-varying parameter observer}

It is obvious that solution of \eqref{eta} is a harmonical signal and we can write down the solution of the equation (\ref{eta}). Let us suppose on the first step that $\omega$ is known 
\begin{equation}
	\eta(t)=a_{1}\sin(\omega \phi(t))+a_{2}\cos(\omega \phi(t)),
	\label{eta_sol}
\end{equation}
where $ a_{1} $ and $a_{2}$ are unknown constant parameters. Let us denote
\begin{equation}
	\chi=\begin{bmatrix}
		\sin (\omega \phi(t))\\ 
		\cos (\omega \phi(t)) \end{bmatrix}  .
	\label{phi}
\end{equation}

Then (\ref{eta_sol}) can be rewritten in the following form:
\begin{equation}
	\eta(t)=a_{1}{\chi }_{1}+a_{2}{\chi }_{2}
\end{equation}
or
\begin{equation}
	\eta (t)=a^\top \chi ,
	\label{eta_matr}
\end{equation}
where $a =\begin{bmatrix}
	{a}_{1}\\ 
	{a }_{2} \end{bmatrix}.$

Substituting (\ref{eta_matr}) into (\ref{sys_delay}), we obtain:
\begin{equation}
	\dot{{x}}_d={A}_d {x}_d+{a}^\top \chi {x}_d+{B}_d {u}_d.
\end{equation}

Let us consider the filter $ \frac{\lambda_2 }{p+\lambda_2 } $. If we apply it to previous equation, we can transform the initial system into linear regression form
\begin{equation}
	\caly ={a }_{1}{\psi }_{1}+{a }_{2}{\psi }_{2},
\end{equation}
where
\begin{equation}
	\caly := \begin{bmatrix} {\caly }_{1} \\ 
		{\caly }_{2}\end{bmatrix} = \frac{\lambda_2 p}{p+\lambda_2 }\begin{bmatrix}{{x}_d}_{1}\\ 
		{{x}_d}_{2}\end{bmatrix} -\frac{\lambda_2 }{p+\lambda_2 }({A}_d\begin{bmatrix}{{x}_d}_{1}\\ 
		{{x}_d}_{2}\end{bmatrix} +{B}_d {u}_d),
\end{equation}
\begin{equation}
	{\psi }_{1}:=\begin{bmatrix}{\psi }_{11}\\ 
		{\psi }_{12}\end{bmatrix} =\frac{\lambda_2 }{p+\lambda_2 }{\chi }_{1}\begin{bmatrix}{{x}_d}_{1}\\ 
		{{x}_d}_{2}\end{bmatrix},
\end{equation}
\begin{equation}
	{\psi }_{2}:=\begin{bmatrix}{\psi }_{21}\\ 
		{\psi }_{22}\end{bmatrix} =\frac{\lambda_2 }{p+\lambda_2 }{\chi }_{2}\begin{bmatrix}{{x}_d}_{1}\\ 
		{{x}_d}_{2}\end{bmatrix}.
\end{equation}

The linear regression model can be rewritten as a system of linear equations
\begin{equation}
	\begin{cases} 
		\label{lin_reg}
		\caly_{1}={a }_{1}{\psi }_{11}+{a }_{2}{\psi }_{21}, \\ 
		{\caly }_{2}={a }_{1}{\psi }_{12}+{a }_{2}{\psi }_{22}.
	\end{cases}
\end{equation}

For estimation of unknown parameters $a_{1}$ and $a_{2}$ we suggest using dynamical regression extention and mixing (DREM) technology \cite{b7,b10, b11} in the form it was developed in \cite{b8}. Also we can use only one equation from \eqref{lin_reg}. Then the observer for unknown parameters $a_{1}$ and $a_{2}$ can be written in the following form 
\begin{subequations}
\label{DREM}
\begin{align}
	\label{doty}
		\dot Y &= - \lambda_3 Y +  \lambda\Psi^\top \caly_1,\\
	\label{dotome}
		\dot \Omega &=- \lambda_3 \Omega +  \lambda\Psi^\top \Psi, \\
	\dot {\hat a} &=-\gamma_2 \Delta (\Delta \hat{\theta}-\calz),
		\label{drem_grad}
	\end{align}
\end{subequations}

with $\lambda _3>0$ and $\gamma_2>0$, $\hat{a}=
\begin{bmatrix}
	\hat{a}_1 \\ \hat{a}_2
\end{bmatrix}$ with the definitions
\begin{subequations}
	\label{gpebodyn1}
	\begin{align}
		\label{Psi}
		\Psi&:=
		\begin{bmatrix}
			\psi_{11} \; \psi_{21}
		\end{bmatrix},\\
		\label{caly}
		\calz &:= \adj\{\Omega\}Y,\\
		\Delta&:=\det\{\Omega\},
		\label{del}
	\end{align}
\end{subequations}

Substituting $ {\hat{a }}_{1}, {\hat{a }}_{2} $ into  (\ref{eta_sol}), we have
\begin{equation}
	\hat{\eta }(t)={\hat{a }}_{1}\sin (\omega \phi(t) )+{\hat{a }}_{2}\cos (\omega  \phi(t)).
\end{equation}

If we use (\ref{eta_theta}), we can obtain the unknown time-varying parameter $ \theta(t) $ by substitution of $ {\hat{a }}_{1}, {\hat{a }}_{2} $ and $\hat{\omega}$
\begin{align}
	\hat{\theta }(t)&=\hat{\eta }(t+d)={\hat{a }}_{1}\sin (\hat{\omega }t)+{\hat{a }}_{2}\cos (\hat{\omega }t). 
	\end{align}

Now we can find estimates of the state vector of \eqref{sys}, \eqref{out} using algorithm with generalized parameter estimation-based observers (GPEBO) technique \cite{b18} for delayed systems  from \cite{b2}. If we use the estimate of $\hat{\theta}(t)$ then we have the following identification algorithm.
\begin{proposition}
	Consider a dynamical system
	\begin{subequations}
		\label{gpebodyn}
		\begin{align}
			\label{dotxi}
			\dot{\xi}(t)&=A(t)\xi(t)+\hat{\theta}(t)x(t)+B(t)u(t),\\
			\label{dotphi}
			\dot{\Phi}_A(t)&=\left( A(t)+\hat{\theta}(t) I_n \right) \Phi_A(t),\;\Phi_A(0)=I_n,
			%\label{doty}
			%\dot Y(t) &= - \lambda Y(t) +  \lambda\Psi^\top (t)z(t)\\
			%\label{dotome}
			%\dot \Omega(t) &=- \lambda \Omega(t) +  \lambda %\Psi^\top(t)\Psi(t)
		\end{align}
	\end{subequations}
	and the gradient parameter estimator
	\begin{align}
		\label{dothatthe}
		\dot {\hat e}(t) &=-\gamma_3 	\calp(t) [	\calp \hat{e}(t)-\calr(t)],
	\end{align}
	with $\gamma_3>0$.
	
	Let us define $\calr(t)$ and $\calp(t)$ as
	\begin{subequations}
		\label{gpebodyn1}
		\begin{align}
			\label{calz}
			\calr(t) &:= \adj\{\Phi_A(\phi(t))\}q(t),\\
			\calp(t)&:=\det\{\Phi_A(\phi(t))\},
			\label{del}
		\end{align}
	\end{subequations}
	where $\adj\{\cdot\}$ is the adjugate matrix. 
	
	Then we define the state estimate as
	\begin{subequations}
	\begin{align}
		\label{ftcgpebo}
		\hat x(t) &= \xi(t) {+} \Phi_A(t) \hat{e}_{FT},\\
		\label{e_ft}
		\hat{e}_{FT} &= {1 \over 1 - w_c(t)}[\hat e(t) - w_c(t) \hat e(0)]
	\end{align}
	\end{subequations}
	with
	\begin{align}
		\dot w(t)  &= -\gamma \Delta^2(t) w(t), \; w(0)=1,
		\label{dotw}
	\end{align}
	and $w_c(t)$ defined via the clipping function
	\begin{equation}
		\label{wc}
		w_c(t) = \left\{ \begin{array}{lcl} w(t) & \;\mbox{if}\; & w(t) \leq 1-\mu, \\ 1-\mu & \;\mbox{if}\; & w(t) > 1-\mu, \end{array} \right.
	\end{equation}
	where $\mu \in (0,1)$ is a designer chosen parameter.

	Then $\hat x(t) = x(t),\;\forall t \geq t_c,$ for some $t_c \in (0,\infty)$
\end{proposition}

\begin{proof}
	
Consider error equation
\begin{align}
	\label{error}
	e(t)=x(t)-\xi(t).
\end{align}
Then, taking into account $\hat{\theta}(t)=\theta(t)$ we will have
\begin{align}
\dot{e}(t)=\left( A(t)- \theta(t) I_n  \right) e(t).
\end{align}
The solution for $e(t)$ can be found in the following form
\begin{align}
	\label{e}
	e(t)=\Phi_A e(0),
\end{align}
where $e(0)=x(0) - \xi(0)$ and $\Phi_A$ is the fundamental matrix defined by equation \eqref{dotphi}. For zero initial conditions in \eqref{dotxi} we have  $e(0)=x(0)$.

After substitution \eqref{e} into \eqref{error} we can write
\begin{align}
	\label{error_1}
	\Phi_A(t) e(0)=x(t)-\xi(t).
\end{align}
In \eqref{error_1} the unknown state vector $x$ is used, but we can implement equation written in the following form
\begin{align}
	%\label{error_1}
	\Phi_A(\phi(t)) e(0)=x(\phi(t))-\xi(\phi(t)),
\end{align}
and we obtain linear regression equation
\begin{align}
	%\label{error_1}
q(t)=	\Phi_A(\phi(t)) e(0),
\end{align}
where $q(t)=x(\phi(t))-\xi(\phi(t))$.
Now we can find vector of initial conditions $(\hat{e}(0))$ using gradient algorithm \eqref{dothatthe}  or $(\hat{e}_{FT}(0))$ using finite time  algorithm \eqref{e_ft}. After estimation of initial condition the state vector can be found by \eqref{ftcgpebo}.
\end{proof}

%\begin{align}
%	\dot{\hat{x}}(t)&=A(t)\hat{x}(t)+\hat{\theta}(t)\hat{x}(t)+B(t)u(t).
%	\end{align}

\section{Simulation Results}

We consider system \eqref{sys}, \eqref{out} with the following parameters \\
$A(t)=
\begin{bmatrix}
	0 \; &1+0.1 \sin(t)\\
	-2 \; &-1+0.5 \cos(2t)
\end{bmatrix}$, $B=
\begin{bmatrix}
	0 \\ 1
\end{bmatrix}$, $C=
\begin{bmatrix}
	1 \; &0\\0 \;&1 
\end{bmatrix}$.
For simulation we used initial conditions $x(0)=
\begin{bmatrix}
	1 \\ 2
\end{bmatrix}$, 
and 
$\theta(0)=
\begin{bmatrix}
	1.732 \\ 3
\end{bmatrix}$,
$\omega=3$, $u=2 \sin (t)$.
We use $\lambda_1=\lambda_2=\lambda_3=10$ as filter parameter. 

\subsection{System without time delay}

\begin{figure}[hbtp]
	\begin{center}
		\includegraphics[width=1 \linewidth]{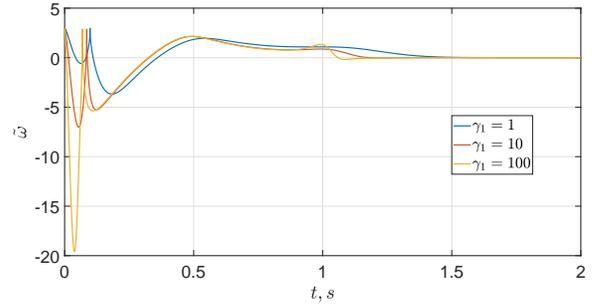}
		\caption{Transients of the error $\tilde{\omega}=\omega-\hat{\omega}$ for different values of adaptation gain $\gamma_1$ (case without delay)}
		\label{err_w_no_del}
	\end{center}
\end{figure}

\begin{figure}[hbtp]
	\begin{center}
		\includegraphics[width=1 \linewidth]{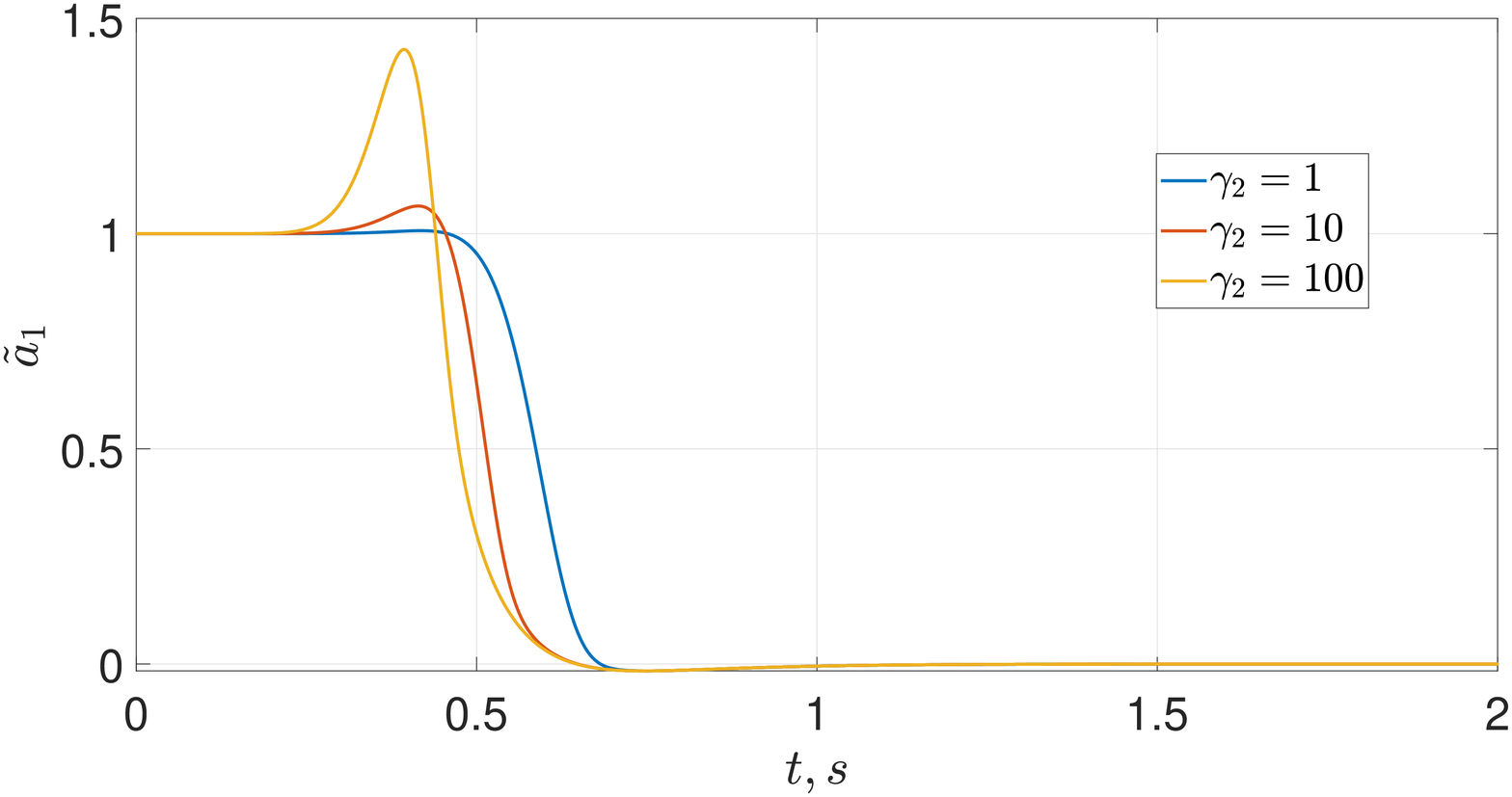}
		\caption{Transients of the error $\tilde{a}_1=a_1-\hat{a}_1$ for case when $\omega$ is known and different values of adaptation gain $\gamma_2$ (case without delay)}
		\label{err_a1_no_del}
	\end{center}
\end{figure}

\begin{figure}[hbtp]
	\begin{center}
		\includegraphics[width=1 \linewidth]{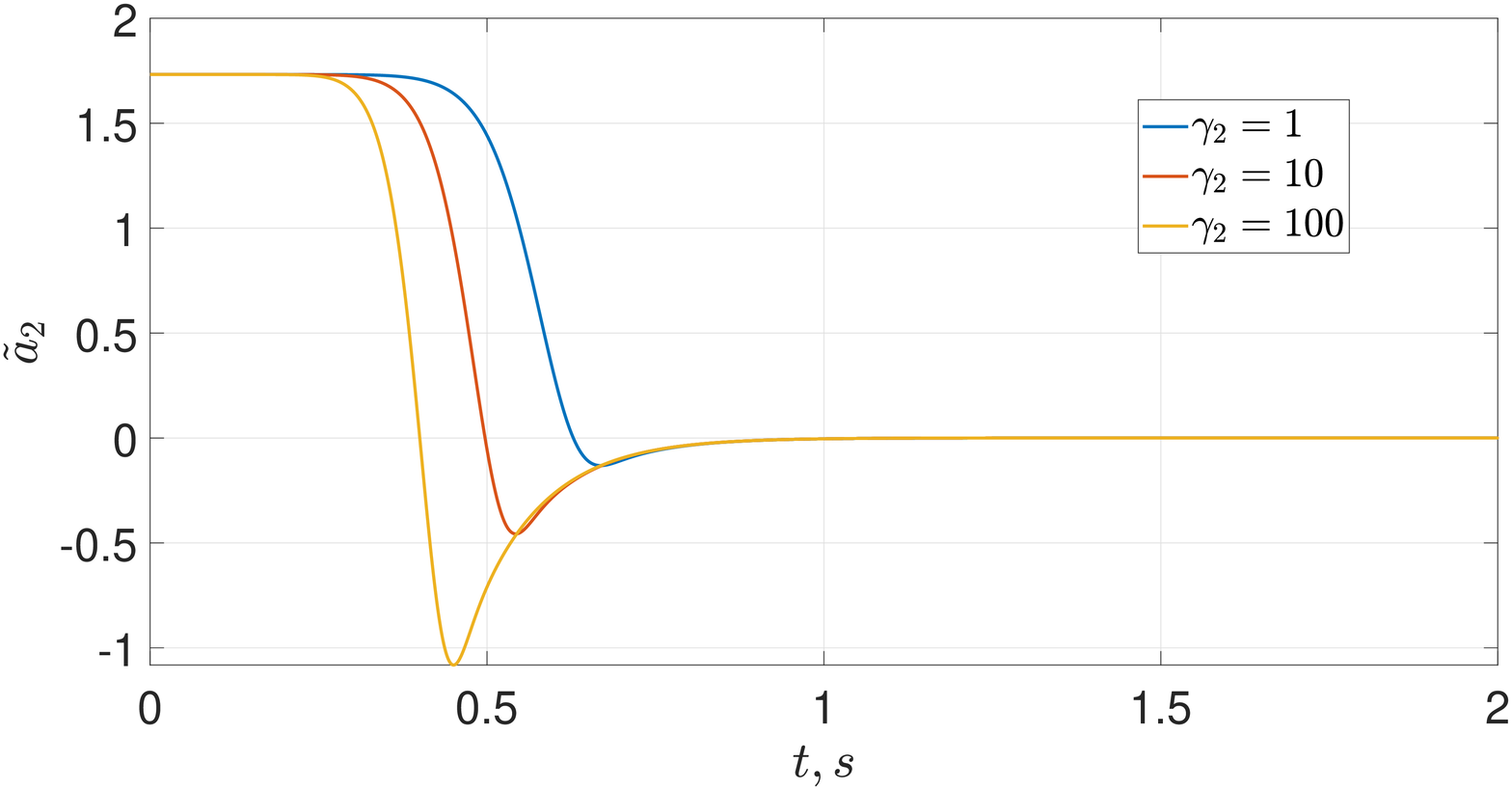}
		\caption{Transients of the error $\tilde{a}_2=a_2-\hat{a}_2$ for case when $\omega$ is known and different values of adaptation gain $\gamma_2$ (case without delay)}
		\label{err_a2_no_del}
	\end{center}
\end{figure}

\begin{figure}[hbtp]
	\begin{center}
		\includegraphics[width=1 \linewidth]{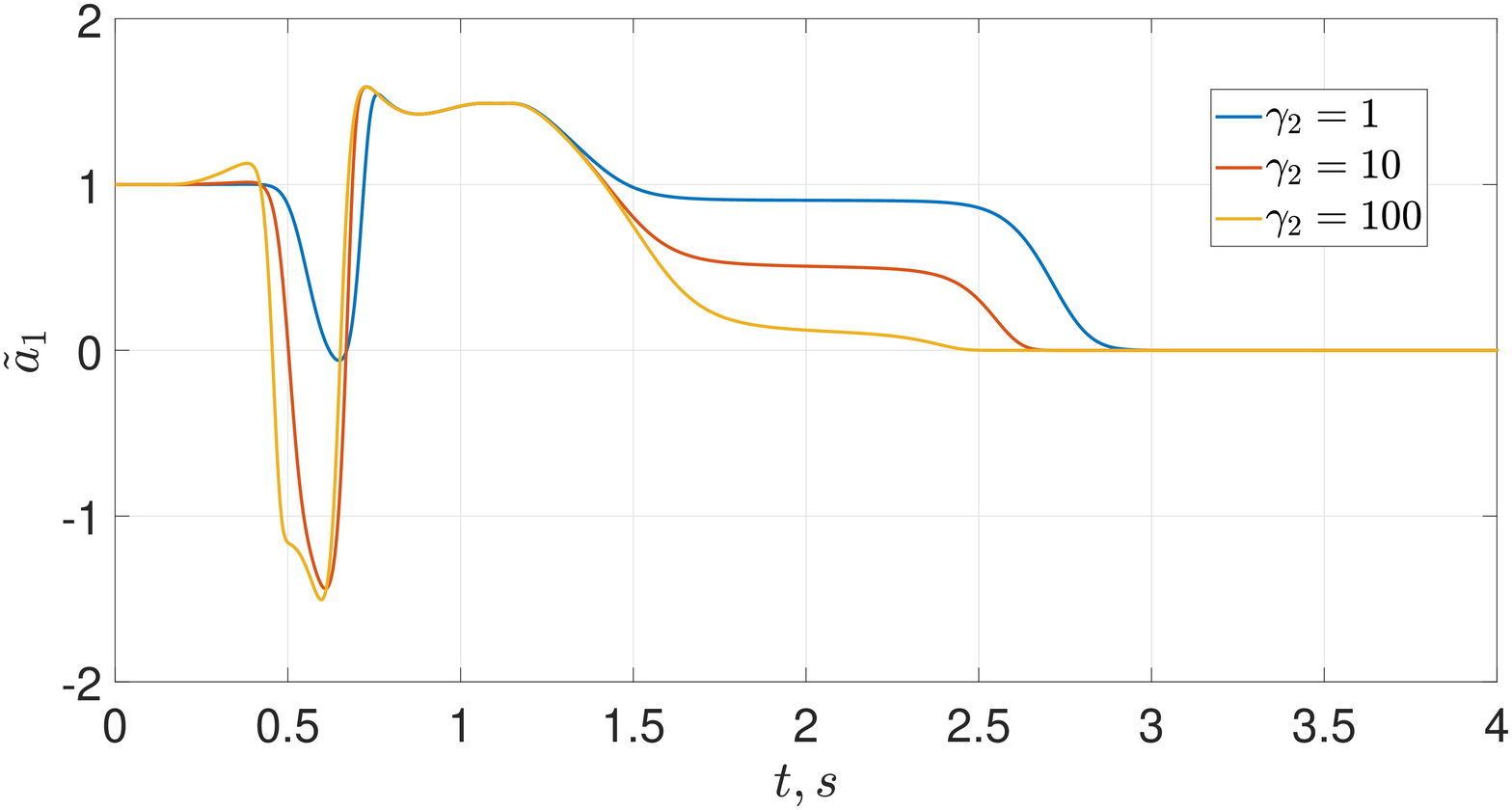}
		\caption{Transients of the error $\tilde{a}_1=a_1-\hat{a}_1$ for case when $\omega$ is unknown and different values of adaptation gain $\gamma_2$ (case without delay)}
		\label{err_a1_1_no_del}
	\end{center}
\end{figure}

\begin{figure}[hbtp]
	\begin{center}
		\includegraphics[width=1 \linewidth]{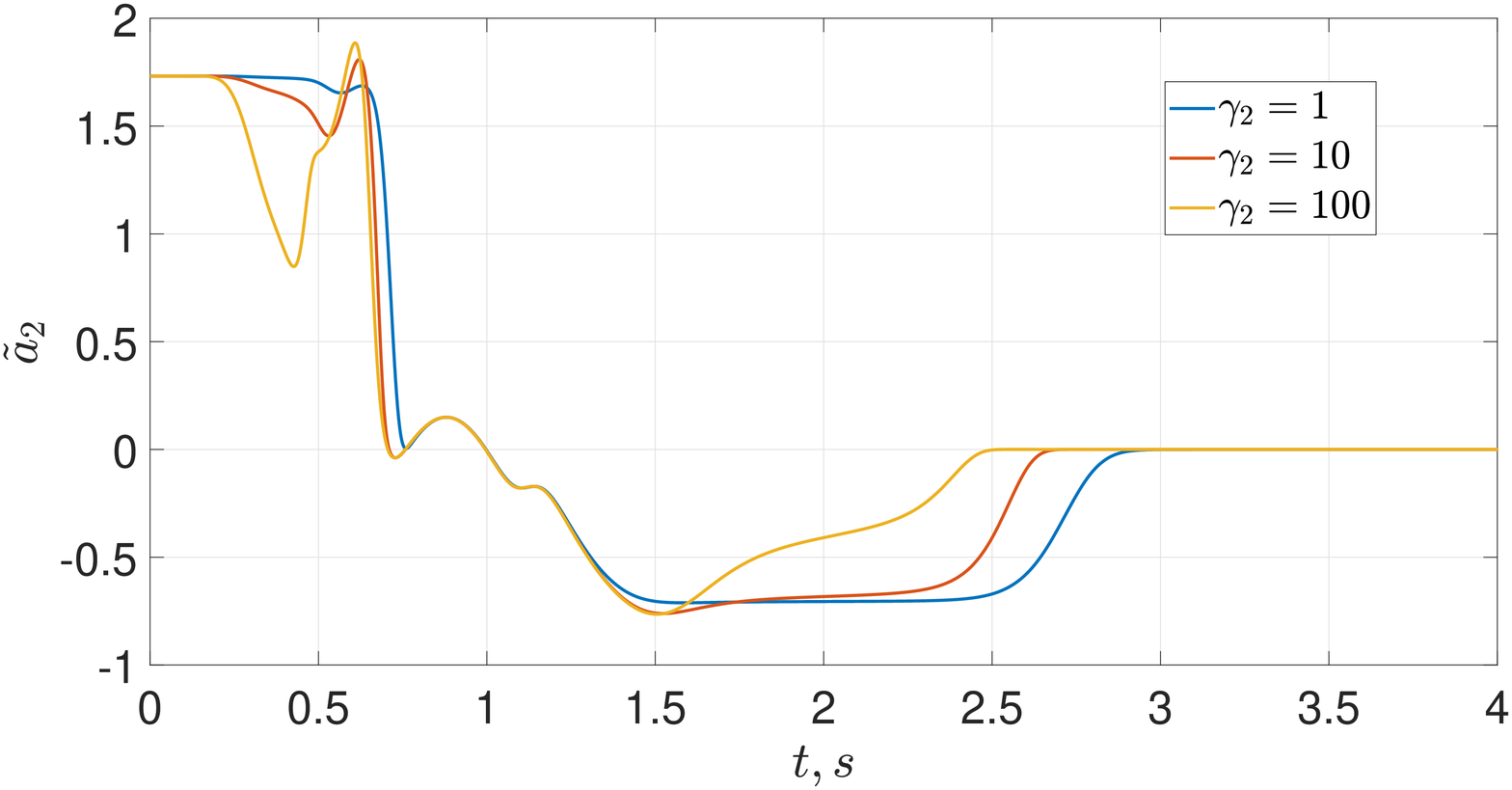}
		\caption{Transients of the error $\tilde{a}_2=a_2-\hat{a}_2$ for case when $\omega$ is unknown and different values of adaptation gain $\gamma_2$ (case without delay)}
		\label{err_a2_1_no_del}
	\end{center}
\end{figure}

\begin{figure}[hbtp]
	\begin{center}
		\includegraphics[width=1 \linewidth]{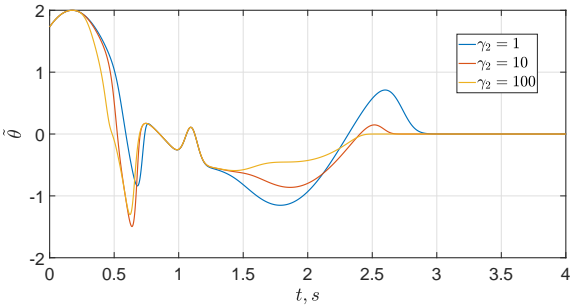}
		\caption{Transients of the error $\tilde{\theta}=\theta-\hat{\theta}$ for different values of adaptation gain $\gamma_2$ (case without delay)}
		\label{err_theta_no_del}
	\end{center}
\end{figure}

Fig.~\ref{err_w_no_del}...Fig.~\ref{err_theta_no_del} demonstrate transients of the proposed algorithms for LTV system \eqref{sys} without delay $(d=0)$. Fig.~\ref{err_w_no_del} demonstrates transients of identification error $\tilde{\theta}=\theta-\hat{\theta}$ for different values of adaptation gain $\gamma_1$. Fig.~\ref{err_a1_no_del}...Fig.~\ref{err_a2_1_no_del} demonstrate transients of identification errors for unknown coefficients $a_1$ and $a_2$ of unknown function $\theta$ \eqref{theta} for different values of adaptation gain $\gamma_2$ for two cases: the first case - we suppose that parameter $\omega$ in \eqref{theta} is known (Fig. \ref{err_a1_no_del} and \ref{err_a2_no_del}) and the second case - we use the estimated value of $\omega$ found by \eqref{hat_omega} and fixed value of adaptation gain $\gamma_1=10$ (Fig. \ref{err_a1_1_no_del} and \ref{err_a2_1_no_del}). Fig.~\ref{err_theta_no_del} demonstrates transients of identification error $\tilde{\theta}=\theta-\hat{\theta}$ for different values of adaptation gain $\gamma_2$ and fixed value of adaptation gain $\gamma_1=10$.

\subsection{System with time delay}

\begin{figure}[hbtp]
	\begin{center}
		\includegraphics[width=1 \linewidth]{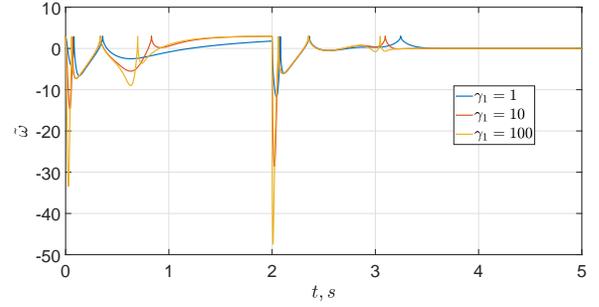}
		\caption{Transients of the error $\tilde{\omega}=\omega-\hat{\omega}$ for different values of adaptation gain $\gamma_1$ (case with delay)}
		\label{err_w_del}
	\end{center}
\end{figure}

\begin{figure}[hbtp]
	\begin{center}
		\includegraphics[width=1 \linewidth]{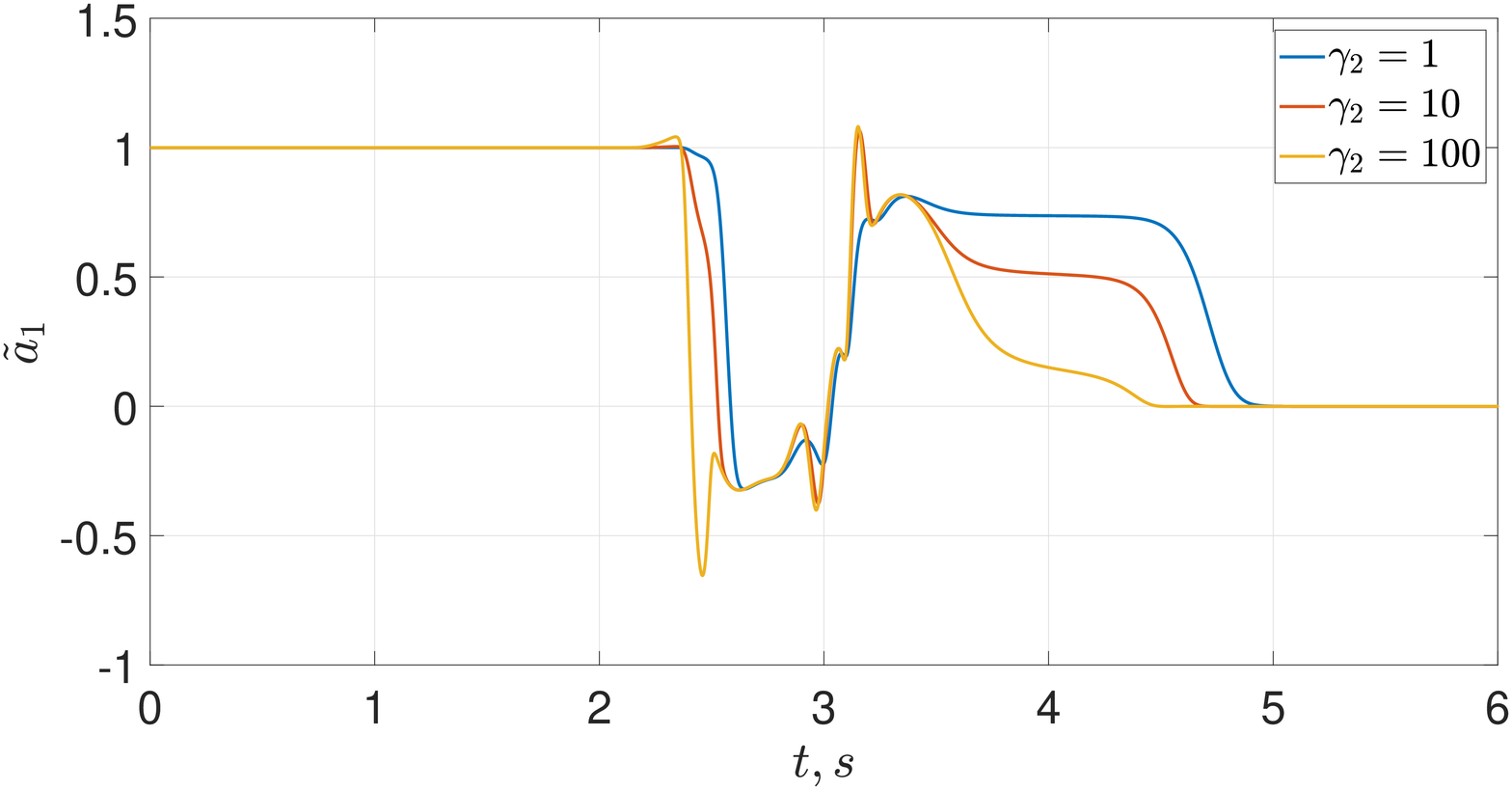}
		\caption{Transients of the error $\tilde{a}_1=a_1-\hat{a}_1$ for case when $\omega$ is unknown and different values of adaptation gain $\gamma_2$ (case with delay)}
		\label{err_a1_del}
	\end{center}
\end{figure}

\begin{figure}[hbtp]
	\begin{center}
		\includegraphics[width=1 \linewidth]{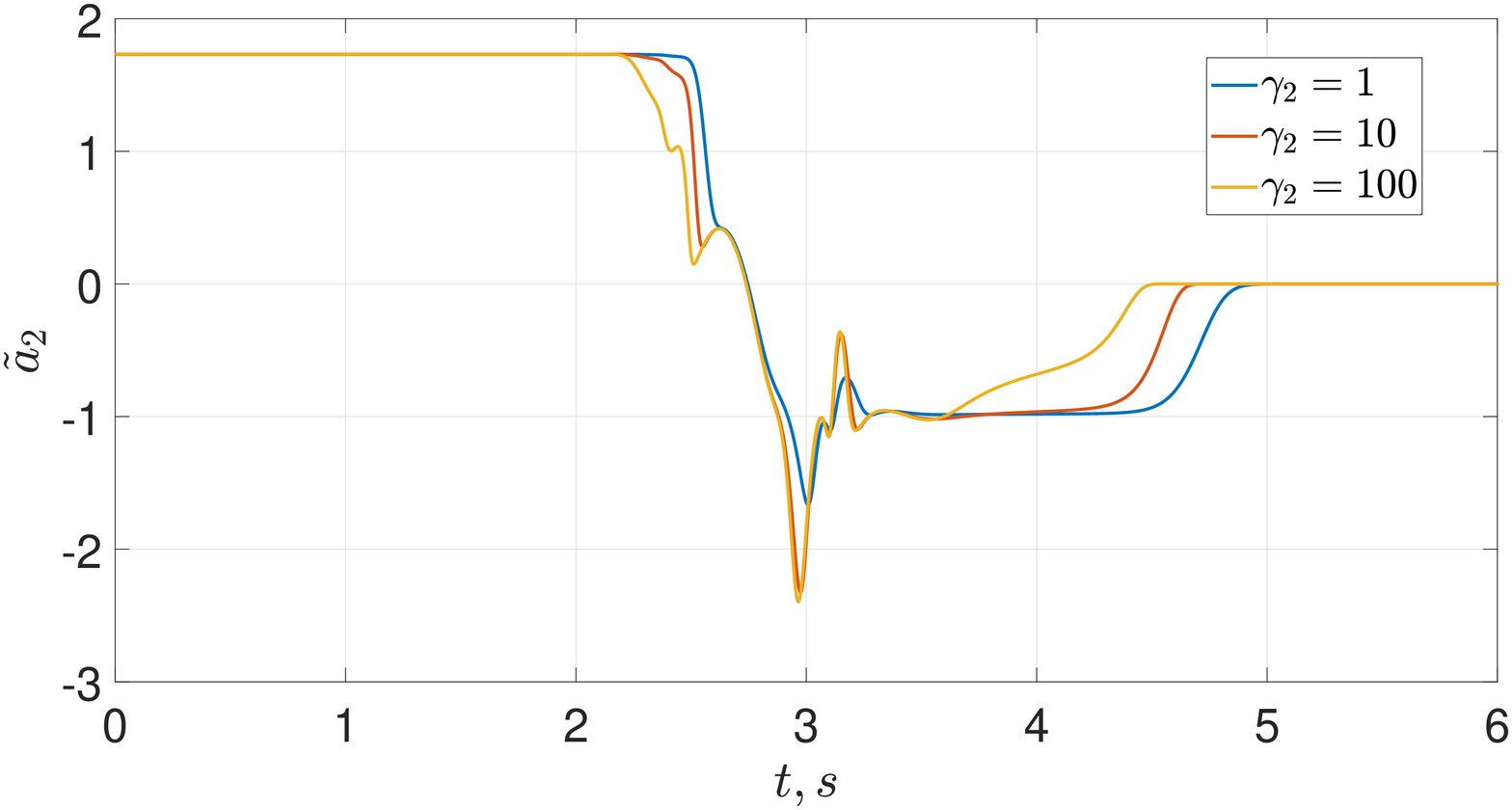}
		\caption{Transients of the error $\tilde{a}_2=a_2-\hat{a}_2$ for case when $\omega$ is unknown and different values of adaptation gain $\gamma_2$ (case with delay)}
		\label{err_a2_del}
	\end{center}
\end{figure}

\begin{figure}[hbtp]
	\begin{center}
		\includegraphics[width=1 \linewidth]{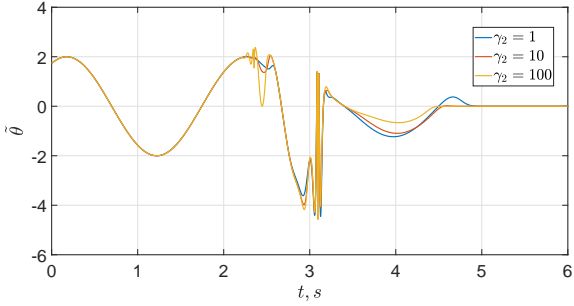}
		\caption{Transients of the error $\tilde{\theta}=\theta-\hat{\theta}$ for different values of adaptation gain $\gamma_2$ (case with delay)}
		\label{err_theta_del}
	\end{center}
\end{figure}

For simulation we used fixed value of time delay $d=2$. Fig.~\ref{err_w_del}...Fig.~\ref{err_theta_del} demonstrate transients of the proposed algorithms. Fig.~\ref{err_w_del} demonstrates transients of identification error $\tilde{\theta}=\theta-\hat{\theta}$ for different values of adaptation gain $\gamma_1$. Fig.~\ref{err_a1_del}...Fig.~\ref{err_a2_del} demonstrate transients of identification errors for unknown coefficients $a_1$ and $a_2$ of unknown function $\theta$ \eqref{theta} when we use the estimated value of $\omega$ obtaned by \eqref{hat_omega} and fixed value of adaptation gain $\gamma_1=10$. Fig.~\ref{err_theta_del} demonstrates transients of identification error $\tilde{\theta}=\theta-\hat{\theta}$ for different values of adaptation gain $\gamma_2$ and fixed value of adaptation gain $\gamma_1=10$. %Fig.~\ref{err_x} demonstrates error transiens for estimates of the state vector of LTV system \eqref{sys}, \eqref{out} with fixed values of adaptation gains $\gamma_1=10$ and $\gamma_2=100$.

%\begin{figure}[hbtp]
%	\begin{center}
%		\includegraphics[width=1 \linewidth]{figures/err_x.eps}
%		\caption{Transients of the error $\tilde{x}_i=x_i-\hat{x}_i$, $i=1,2$ (case with delay)}
%		\label{err_x}
%	\end{center}
%\end{figure}

\subsection{System without time delay. State vector estimation.}
As we have previously estimated the unknown time-varying parameter $\theta(t)$, we start estimation of the state vector by using GPEBO technics. For simulations we used adaptation gain $\gamma_3=$ in \eqref{dothatthe} and $\mu=0.01$ in \eqref{wc}.

\begin{figure}[hbtp]
	\begin{center}
		\includegraphics[width=1 \linewidth]{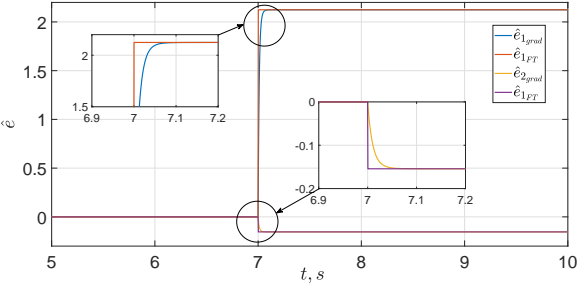}
		\caption{Transients of initial conditions estimations for gradient and finite time observers (case with delay)}
		\label{hat_e}
	\end{center}
\end{figure}

If we look at the Fig. \ref{err_theta_del}, we can see that  $\theta(t)$ was estimated within 5 second and then state vector observer is switched on (see Fig. \ref{hat_e}, which demonstrates transients for estimations of initial conditions).

\begin{figure}[hbtp]
	\begin{center}
		\includegraphics[width=1 \linewidth]{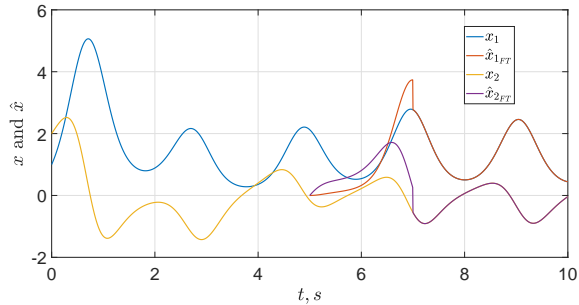}
		\caption{Transients of the state vector estimations with finite time observer (case with delay)}
		\label{x}
	\end{center}
\end{figure}

\section{Concluding Remarks} 

We have presented the adaptive observer for LTV system \eqref{sys}, \eqref{out} with partially unknown state matrix and  delayed measurements. The proposed observer performance consists of two steps. On the first step we estimate the unknown time-varying parameter. On the second step we use GPEBO technique to estimate unknown state vector. The simulation results demonstrate efficiency of the proposed algorithms. The future development of current research can be case of the system which state is not measured and only output is measured with delay and $y=C^\top x(\phi(t))$, where $C$ is a vector.

\vspace{12pt}


\begin{thebibliography}{00}
	%1
	\bibitem{b4} R. Sanx, P. Garcia, M. Krstic, "Observation and stabilization of LTV systems with time-varying measurements delay," Automatica 103, pp.573-579, 2019.
	%2
	\bibitem{b2} A. Bobtsov, N. Nikolaev, R. Ortega, D. Efimov, "State observation of LTV systems with delayed measurements: A parameter estimation-based approach with fixed convergence time," Automatica 131, 109674, 2021.
	%3
	\bibitem{b8} A. Bobtsov, N. Nikolaev, R. Ortega, D. Efimov, "State observation of affine-in-the-states time-varying systems with unknown parameters and delayed measurements," Third IFAC Conference on Modelling, Identification and Control of Nonlinear Systems, pp. 124–129, 2021.
	%4
	\bibitem{b9} A. Bobtsov, N. Nikolaev, R. Ortega, D. Efimov, O. Kozachek, "State Observation of Affine-in-the-States Systems with Unknown Time-Varying Parameters and Output Delay," 14th IFAC International Workshop on Adaptation and Learning in Control and Signal Processing, Casablanca, Morocco, pp.343 - 348, 2022.
	%5
	\bibitem{b1} E. Fridman Introduction to time-delay systems: analysis and control. – Springer, 2014.
	%6
	\bibitem{b3} J. G. Rueda-Escobedo, R. Ushirobira, D. Efimov, J. A. Moreno, "Gramian-based uniform convergent observer for stable LTV systems with delayed measurements," International Journal of Control 93 (2), pp.226-237, 2020.
	%7
	\bibitem{b12} R. A. Biroon, Z. Abdollahi and P. Pisu, "Measurement Unknown Delay Estimation in Cyber-Physical Systems: a PDE approach," 2021 American Control Conference (ACC), 2021, pp. 4619-4624, doi: 10.23919/ACC50511.2021.9483246.
	%8
	\bibitem{b16} B. Guo and Z. Mei, "Output Feedback Stabilization for a Class of First-Order Equation Setting of Collocated Well-Posed Linear Systems With Time Delay in Observation," in IEEE Transactions on Automatic Control, vol. 65, no. 6, pp. 2612-2618, June 2020, doi: 10.1109/TAC.2019.2941431.
	%9
	\bibitem{b13} R. Ramjug-Ballgobin, K. Busawon, R. T. F. Ah King and H. C. S. Rughooputh, "PI Observer-based Control for Biomass Regulation under Measurement Delays," 2021 11th IEEE International Conference on Control System, Computing and Engineering (ICCSCE), 2021, pp. 71-74, doi: 10.1109/ICCSCE52189.2021.9530946.
	%10
	\bibitem{b14} R. Ramjug-Ballgobin, K. Busawon, R. T. F. Ah King and H. C. S. Rughooputh, "Observer-based Control for Biomass Regulation under Discrete Measurements," 2020 10th IEEE International Conference on Control System, Computing and Engineering (ICCSCE), 2020, pp. 19-22, doi: 10.1109/ICCSCE50387.2020.9204945.
	%11
	\bibitem{b15} R. P. Desai and N. S. Manjarekar, "Pitch Channel Tracking Control of an Autonomous Underwater Vehicle with Delayed Output," 2021 IEEE 18th India Council International Conference (INDICON), 2021, pp. 1-6, doi: 10.1109/INDICON52576.2021.9691723.
	%12
	\bibitem{b17} N. Sehli, K. Ibn Taarit, Z. Wang, T. Raïssi and M. Ksouri, "Joint interval state and actuator fault estimation for linear discrete-time delayed systems," 2021 International Conference on Control, Automation and Diagnosis (ICCAD), 2021, pp. 1-6, doi: 10.1109/ICCAD52417.2021.9638746.
	%13
	\bibitem{b5} L. Ljung System identification: theory for the users, Prentice Hall, New Jersey, 1999.
	%14
	\bibitem{b6} S. Sastry, M. Bodson, J. F. Bartram, Adaptive control: stability, convergence, and robustness, 1990.
	%15
	\bibitem{b7} S. Aranovskiy, A. Bobtsov, R. Ortega and A. Pyrkin, "Performance Enhancement of Parameter Estimators via Dynamic Regressor Extension and Mixing," IEEE Transactions on Automatic Control, vol. 62, no. 7, pp. 3546-3550, July 2017, doi: 10.1109/TAC.2016.2614889.
	%16
	\bibitem{b10} S. Aranovskiy, A. Bobtsov, A. Pyrkin, R. Ortega, A. Chaillet," Flux and Position Observer of Permanent Magnet Synchronous Motors with Relaxed Persistency of Excitation Conditions," IFAC-PapersOnLine, 48(11), pp. 301–306, 2015.
	%17
	\bibitem{b11} R. Ortega, S. Aranovskiy, A. Pyrkin, A. Astolfi, A. Bobtsov, "New Results on Parameter Estimation via Dynamic Regressor Extension and Mixing: Continuous and Discrete-Time Cases," IEEE Transactions on Automatic Control,  66(5), pp. 2265–2272, 2021.
	%18	
	\bibitem{b18} R. Ortega, A. Bobtsov, N. Nikolaev, J. Schiffer, D. Dochain. "Generalized parameter estimation-based observers: Application to power systems and chemical–biological reactors," Automatica 129, 109635, 2021.
	
\end{thebibliography}
\end{document}